\documentclass[12pt]{amsart}
\usepackage{amssymb}
\begin{document}
%%%%%%%%%%%%%%%%%%%% Text italic %%%%%%%%%%%%%%%%%%%%%%%%%%%%
\theoremstyle{plain}
\newtheorem{thm}{Theorem}[section]

\newtheorem{theorem}[thm]{Theorem}
\newtheorem{lemma}[thm]{Lemma}
\newtheorem{corollary}[thm]{Corollary}
\newtheorem{proposition}[thm]{Proposition}
\newtheorem{addendum}[thm]{Addendum}
\newtheorem{variant}[thm]{Variant}

\newtheorem{mythm}{Theorem}
\newtheorem{myprop}[mythm]{Proposition}
\newtheorem{mylem}[mythm]{Lemma}
\newtheorem{mycor}[mythm]{Corollary}

\newtheorem*{mythm1}{Theorem}
\newtheorem*{myprop1}{Proposition}
\newtheorem*{mylem1}{Lemma}
\newtheorem*{mycor1}{Corollary}
\newtheorem*{myexe1}{Exercise}

\newtheorem{mythma}{\bf{Theorem}}[section]
\newtheorem{mypropa}[mythma]{\bf{Proposition}}
\newtheorem{mylema}[mythma]{\bf{Lemma}}
\newtheorem{mycora}[mythma]{Corollary}
\newtheorem{myexea}[mythma]{\bf{Exercise}}

\newtheorem{mythmb}{Theorem}[subsection]
\newtheorem{theoremnum}[mythmb]{Theorem}
\newtheorem{mycorb}[mythmb]{Corollary}
\newtheorem{mypropb}[mythmb]{Proposition}
\newtheorem{myexeb}[mythmb]{\bf{Exercise}}

%%%%%%%%%%%%%%%%%%%% Text roman %%%%%%%%%%%%%%%%%%%%%%%%%%%%%
\theoremstyle{definition}
\newtheorem{construction}[thm]{Construction}
\newtheorem{conjecture}[thm]{Conjecture}
\newtheorem{notations}[thm]{Notations}
\newtheorem{notation}[thm]{Notation}
\newtheorem{question}[thm]{Question}
\newtheorem{problem}[thm]{Problem}
\newtheorem{remark}[thm]{Remark}
\newtheorem{remarks}[thm]{Remarks}
\newtheorem{definition}[thm]{Definition}
\newtheorem{claim}[thm]{Claim}
\newtheorem{assumption}[thm]{Assumption}
\newtheorem{assumptions}[thm]{Assumptions}
\newtheorem{properties}[thm]{Properties}
\newtheorem{example}[thm]{Example}
\newtheorem{examples}[thm]{Examples}
\newtheorem{exercise}[thm]{Exercise}
\newtheorem{comments}[thm]{Comments}
\newtheorem{blank}[thm]{}
\newtheorem{fact}[thm]{Fact}
\newtheorem{defn-thm}[thm]{Definition-Theorem}
\newtheorem{criterion}[thm]{Criterion}

\newtheorem*{myrem}{Remark}
\newtheorem*{myrems}{Remarks}
\newtheorem*{mydef}{Definition}
\newtheorem*{myexam}{Example}
\newtheorem*{myexams}{Examples}
\newtheorem*{myexe}{Exercise}
\newtheorem*{myhint}{Hint}
\newtheorem*{mynote}{Notation}
\newtheorem*{mynotes}{Notations}
\newtheorem*{myack}{Acknowledgments}         %\renewcommand{\theack}{}
\newtheorem{myacknum}[thm]{Acknowledgments}

\numberwithin{equation}{thm}

\newtheorem{mydef1}{Definition}
\newtheorem{mydefa}[thm]{Definition}

%%%%%%%% Diagram macros, etc. %%%%%%%%%%%%%%%%%%%%%%%%%%%%%%%
\catcode`\@=11
% General macros
\def\opn#1#2{\def#1{\mathop{\kern0pt\fam0#2}\nolimits}}
\def\bold#1{{\bf #1}}%
\def\underrightarrow{\mathpalette\underrightarrow@}
\def\underrightarrow@#1#2{\vtop{\ialign{$##$\cr
 \hfil#1#2\hfil\cr\noalign{\nointerlineskip}%
 #1{-}\mkern-6mu\cleaders\hbox{$#1\mkern-2mu{-}\mkern-2mu$}\hfill
 \mkern-6mu{\to}\cr}}}
\let\underarrow\underrightarrow
\def\underleftarrow{\mathpalette\underleftarrow@}
\def\underleftarrow@#1#2{\vtop{\ialign{$##$\cr
 \hfil#1#2\hfil\cr\noalign{\nointerlineskip}#1{\leftarrow}\mkern-6mu
 \cleaders\hbox{$#1\mkern-2mu{-}\mkern-2mu$}\hfill
 \mkern-6mu{-}\cr}}}
% Rectangular Commutative diagrams
\let\amp@rs@nd@\relax
\newdimen\ex@
\ex@.2326ex
\newdimen\bigaw@
\newdimen\minaw@
\minaw@16.08739\ex@
\newdimen\minCDaw@
\minCDaw@2.5pc
\newif\ifCD@
\def\minCDarrowwidth#1{\minCDaw@#1}
\newenvironment{CD}{\@CD}{\@endCD}
\def\@CD{\def\A##1A##2A{\llap{$\vcenter{\hbox
 {$\scriptstyle##1$}}$}\Big\uparrow\rlap{$\vcenter{\hbox{%
$\scriptstyle##2$}}$}&&}%
\def\V##1V##2V{\llap{$\vcenter{\hbox
 {$\scriptstyle##1$}}$}\Big\downarrow\rlap{$\vcenter{\hbox{%
$\scriptstyle##2$}}$}&&}%
\def\={&\hskip.5em\mathrel
 {\vbox{\hrule width\minCDaw@\vskip3\ex@\hrule width
 \minCDaw@}}\hskip.5em&}%
\def\verteq{\Big\Vert&&}%
\def\noarr{&&}%
\def\vspace##1{\noalign{\vskip##1\relax}}\relax\iffalse{%
\fi\let\amp@rs@nd@&\iffalse}\fi
 \CD@true\vcenter\bgroup\relax\iffalse{%
\fi\let\\=\cr\iffalse}\fi\tabskip\z@skip\baselineskip20\ex@
 \lineskip3\ex@\lineskiplimit3\ex@\halign\bgroup
 &\hfill$\m@th##$\hfill\cr}
\def\@endCD{\cr\egroup\egroup}
% Horizontal arrows with "sliding" length
\def\>#1>#2>{\amp@rs@nd@\setbox\z@\hbox{$\scriptstyle
 \;{#1}\;\;$}\setbox\@ne\hbox{$\scriptstyle\;{#2}\;\;$}\setbox\tw@
 \hbox{$#2$}\ifCD@
 \global\bigaw@\minCDaw@\else\global\bigaw@\minaw@\fi
 \ifdim\wd\z@>\bigaw@\global\bigaw@\wd\z@\fi
 \ifdim\wd\@ne>\bigaw@\global\bigaw@\wd\@ne\fi
 \ifCD@\hskip.5em\fi
 \ifdim\wd\tw@>\z@
 \mathrel{\mathop{\hbox to\bigaw@{\rightarrowfill}}\limits^{#1}_{#2}}\else
 \mathrel{\mathop{\hbox to\bigaw@{\rightarrowfill}}\limits^{#1}}\fi
 \ifCD@\hskip.5em\fi\amp@rs@nd@}
\def\<#1<#2<{\amp@rs@nd@\setbox\z@\hbox{$\scriptstyle
 \;\;{#1}\;$}\setbox\@ne\hbox{$\scriptstyle\;\;{#2}\;$}\setbox\tw@
 \hbox{$#2$}\ifCD@
 \global\bigaw@\minCDaw@\else\global\bigaw@\minaw@\fi
 \ifdim\wd\z@>\bigaw@\global\bigaw@\wd\z@\fi
 \ifdim\wd\@ne>\bigaw@\global\bigaw@\wd\@ne\fi
 \ifCD@\hskip.5em\fi
 \ifdim\wd\tw@>\z@
 \mathrel{\mathop{\hbox to\bigaw@{\leftarrowfill}}\limits^{#1}_{#2}}\else
 \mathrel{\mathop{\hbox to\bigaw@{\leftarrowfill}}\limits^{#1}}\fi
 \ifCD@\hskip.5em\fi\amp@rs@nd@}
% Rectangular commutative diagrams with diagonal arows
\newenvironment{CDS}{\@CDS}{\@endCDS}
\def\@CDS{\def\A##1A##2A{\llap{$\vcenter{\hbox
 {$\scriptstyle##1$}}$}\Big\uparrow\rlap{$\vcenter{\hbox{%
$\scriptstyle##2$}}$}&}%
\def\V##1V##2V{\llap{$\vcenter{\hbox
 {$\scriptstyle##1$}}$}\Big\downarrow\rlap{$\vcenter{\hbox{%
$\scriptstyle##2$}}$}&}%
\def\={&\hskip.5em\mathrel
 {\vbox{\hrule width\minCDaw@\vskip3\ex@\hrule width
 \minCDaw@}}\hskip.5em&}
\def\verteq{\Big\Vert&}
\def\novarr{&}
\def\noharr{&&}
\def\SE##1E##2E{\slantedarrow(0,18)(4,-3){##1}{##2}&}
\def\SW##1W##2W{\slantedarrow(24,18)(-4,-3){##1}{##2}&}
\def\NE##1E##2E{\slantedarrow(0,0)(4,3){##1}{##2}&}
\def\NW##1W##2W{\slantedarrow(24,0)(-4,3){##1}{##2}&}
\def\slantedarrow(##1)(##2)##3##4{%
\thinlines\unitlength1pt\lower 6.5pt\hbox{\begin{picture}(24,18)%
\put(##1){\vector(##2){24}}%
\put(0,8){$\scriptstyle##3$}%
\put(20,8){$\scriptstyle##4$}%
\end{picture}}}
\def\vspace##1{\noalign{\vskip##1\relax}}\relax\iffalse{%
\fi\let\amp@rs@nd@&\iffalse}\fi
 \CD@true\vcenter\bgroup\relax\iffalse{%
\fi\let\\=\cr\iffalse}\fi\tabskip\z@skip\baselineskip20\ex@
 \lineskip3\ex@\lineskiplimit3\ex@\halign\bgroup
 &\hfill$\m@th##$\hfill\cr}
\def\@endCDS{\cr\egroup\egroup}
% Triangular commutative diagrams
\newdimen\TriCDarrw@
\newif\ifTriV@
\newenvironment{TriCDV}{\@TriCDV}{\@endTriCD}
\newenvironment{TriCDA}{\@TriCDA}{\@endTriCD}
\def\@TriCDV{\TriV@true\def\TriCDpos@{6}\@TriCD}
\def\@TriCDA{\TriV@false\def\TriCDpos@{10}\@TriCD}
\def\@TriCD#1#2#3#4#5#6{%
\setbox0\hbox{$\ifTriV@#6\else#1\fi$} \TriCDarrw@=\wd0
\advance\TriCDarrw@ 24pt \advance\TriCDarrw@ -1em
\def\SE##1E##2E{\slantedarrow(0,18)(2,-3){##1}{##2}&}
\def\SW##1W##2W{\slantedarrow(12,18)(-2,-3){##1}{##2}&}
\def\NE##1E##2E{\slantedarrow(0,0)(2,3){##1}{##2}&}
\def\NW##1W##2W{\slantedarrow(12,0)(-2,3){##1}{##2}&}
\def\slantedarrow(##1)(##2)##3##4{\thinlines\unitlength1pt
\lower 6.5pt\hbox{\begin{picture}(12,18)%
\put(##1){\vector(##2){12}}%
\put(-4,\TriCDpos@){$\scriptstyle##3$}%
\put(12,\TriCDpos@){$\scriptstyle##4$}%
\end{picture}}}
\def\={\mathrel {\vbox{\hrule
   width\TriCDarrw@\vskip3\ex@\hrule width
   \TriCDarrw@}}}
\def\>##1>>{\setbox\z@\hbox{$\scriptstyle
 \;{##1}\;\;$}\global\bigaw@\TriCDarrw@
 \ifdim\wd\z@>\bigaw@\global\bigaw@\wd\z@\fi
 \hskip.5em
 \mathrel{\mathop{\hbox to \TriCDarrw@
{\rightarrowfill}}\limits^{##1}}
 \hskip.5em}
\def\<##1<<{\setbox\z@\hbox{$\scriptstyle
 \;{##1}\;\;$}\global\bigaw@\TriCDarrw@
 \ifdim\wd\z@>\bigaw@\global\bigaw@\wd\z@\fi
 \mathrel{\mathop{\hbox to\bigaw@{\leftarrowfill}}\limits^{##1}}
 }
 \CD@true\vcenter\bgroup\relax\iffalse{\fi\let\\=\cr\iffalse}\fi
 \tabskip\z@skip\baselineskip20\ex@
 \lineskip3\ex@\lineskiplimit3\ex@
 \ifTriV@
 \halign\bgroup
 &\hfill$\m@th##$\hfill\cr
#1&\multispan3\hfill$#2$\hfill&#3\\
&#4&#5\\
&&#6\cr\egroup%
\else
 \halign\bgroup
 &\hfill$\m@th##$\hfill\cr
&&#1\\%
&#2&#3\\
#4&\multispan3\hfill$#5$\hfill&#6\cr\egroup \fi}
\def\@endTriCD{\egroup}

\newcounter{Myenumi}
\newenvironment{myenumi}%
{\begin{list}{}{\usecounter{Myenumi}%
\renewcommand{\makelabel}{\arabic{Myenumi}.}%
\settowidth{\leftmargin}{2.n}\settowidth{\labelwidth}{2.n}%
\setlength{\labelsep}{0pt}}}{\end{list}}
\newcounter{Myenumii}
\newenvironment{myenumii}%
{\begin{list}{}{\usecounter{Myenumii}%
\renewcommand{\makelabel}{\alph{Myenumii})}%
\settowidth{\leftmargin}{a)n}\settowidth{\labelwidth}{a)n}%
\setlength{\labelsep}{0pt}}}{\end{list}}
\newcounter{Myenumiii}
\newenvironment{myenumiii}%
{\begin{list}{}{\usecounter{Myenumiii}%
\renewcommand{\makelabel}{\roman{Myenumiii}.}%
\settowidth{\leftmargin}{iv.n}\settowidth{\labelwidth}{iv.n}%
\setlength{\labelsep}{0pt}}}{\end{list}}

\renewenvironment{quote}{\begin{list}{}%
{\settowidth{\leftmargin}{2.n}\setlength{\rightmargin}{0pt}
\renewcommand{\makelabel}{}}%
\item}%
{\end{list}}

\renewenvironment{itemize}%
{\begin{list}{}{\renewcommand{\makelabel}{$\bullet$}%
\settowidth{\leftmargin}{2.n}\settowidth{\labelwidth}{2.n}%
\setlength{\labelsep}{0pt}}}{\end{list}}

\newsymbol\onto 1310
\def\into{\DOTSB\lhook\joinrel\rightarrow}

%%%%%%%%%%%%%%%  End of diagram macros.  %%%%%%%%%%%%%%%%%%%%%%%%%
\newcommand{\sA}{{\mathcal A}}
\newcommand{\sB}{{\mathcal B}}
\newcommand{\sC}{{\mathcal C}}
\newcommand{\sD}{{\mathcal D}}
\newcommand{\sE}{{\mathcal E}}
\newcommand{\sF}{{\mathcal F}}
\newcommand{\sG}{{\mathcal G}}
\newcommand{\sH}{{\mathcal H}}
\newcommand{\sI}{{\mathcal I}}
\newcommand{\sJ}{{\mathcal J}}
\newcommand{\sK}{{\mathcal K}}
\newcommand{\sL}{{\mathcal L}}
\newcommand{\sM}{{\mathcal M}}
\newcommand{\sN}{{\mathcal N}}
\newcommand{\sO}{{\mathcal O}}
\newcommand{\sP}{{\mathcal P}}
\newcommand{\sQ}{{\mathcal Q}}
\newcommand{\sR}{{\mathcal R}}
\newcommand{\sS}{{\mathcal S}}
\newcommand{\sT}{{\mathcal T}}
\newcommand{\sU}{{\mathcal U}}
\newcommand{\sV}{{\mathcal V}}
\newcommand{\sW}{{\mathcal W}}
\newcommand{\sX}{{\mathcal X}}
\newcommand{\sY}{{\mathcal Y}}
\newcommand{\sZ}{{\mathcal Z}}
%mathfrak
\newcommand{\ssA}{{\mathfrak A}}
\newcommand{\ssB}{{\mathfrak B}}
\newcommand{\ssC}{{\mathfrak C}}
\newcommand{\ssD}{{\mathfrak D}}
\newcommand{\ssE}{{\mathfrak E}}
\newcommand{\ssF}{{\mathfrak F}}
\newcommand{\ssG}{{\mathfrak G}}
\newcommand{\ssH}{{\mathfrak H}}
\newcommand{\ssI}{{\mathfrak I}}
\newcommand{\ssJ}{{\mathfrak J}}
\newcommand{\ssK}{{\mathfrak K}}
\newcommand{\ssL}{{\mathfrak L}}
\newcommand{\ssM}{{\mathfrak M}}
\newcommand{\ssN}{{\mathfrak N}}
\newcommand{\ssO}{{\mathfrak O}}
\newcommand{\ssP}{{\mathfrak P}}
\newcommand{\ssQ}{{\mathfrak Q}}
\newcommand{\ssR}{{\mathfrak R}}
\newcommand{\ssS}{{\mathfrak S}}
\newcommand{\ssT}{{\mathfrak T}}
\newcommand{\ssU}{{\mathfrak U}}
\newcommand{\ssV}{{\mathfrak V}}
\newcommand{\ssW}{{\mathfrak W}}
\newcommand{\ssX}{{\mathfrak X}}
\newcommand{\ssY}{{\mathfrak Y}}
\newcommand{\ssZ}{{\mathfrak Z}}

\newcommand{\A}{{\mathbb A}}
\newcommand{\B}{{\mathbb B}}
\newcommand{\C}{{\mathbb C}}
\newcommand{\D}{{\mathbb D}}
\newcommand{\E}{{\mathbb E}}
\newcommand{\F}{{\mathbb F}}
\newcommand{\G}{{\mathbb G}}
\newcommand{\HH}{{\mathbb H}}
\newcommand{\I}{{\mathbb I}}
\newcommand{\J}{{\mathbb J}}
\newcommand{\K}{{\mathbb K}}
\renewcommand{\L}{{\mathbb L}}
\newcommand{\M}{{\mathbb M}}
\newcommand{\N}{{\mathbb N}}
\renewcommand{\P}{{\mathbb P}}
\newcommand{\Q}{{\mathbb Q}}
\newcommand{\R}{{\mathbb R}}
\newcommand{\bS}{{\mathbb S}}

\newcommand{\T}{{\mathbb T}}
\newcommand{\U}{{\mathbb U}}
\newcommand{\V}{{\mathbb V}}
\newcommand{\W}{{\mathbb W}}
\newcommand{\X}{{\mathbb X}}
\newcommand{\Y}{{\mathbb Y}}
\newcommand{\Z}{{\mathbb Z}}
\newcommand{\id}{{\rm id}}

\newcommand{\rank}{{\rm rank}}
\newcommand{\END}{{\mathbb E}{\rm nd}}
\newcommand{\End}{{\rm End}}
\newcommand{\Hg}{{\rm Hg}}
\newcommand{\tr}{{\rm tr}}
\newcommand{\Tr}{{\rm Trace}}
\newcommand{\Sl}{{\rm Sl}}
\newcommand{\SL}{{\rm SL}}
\newcommand{\Gl}{{\rm Gl}}
\newcommand{\Cor}{{\rm Cor}}
\newcommand{\coker}{{\rm coker}}
\newcommand{\GL}{{\rm GL}}
\newcommand{\MT}{{\rm MT}}
\newcommand{\Hdg}{{\rm Hdg}}
\newcommand{\MTV}{{\rm MTV}}
\newcommand{\SO}{{\rm SO}}
\newcommand{\Sp}{{\rm Sp}}
\newcommand{\SU}{{\rm SU}}
\newcommand{\supp}{{\rm Supp}}
\newcommand{\Hom}{{\rm Hom}}
\newcommand{\Ker}{{\rm Ker}}
\newcommand{\Lie}{{\rm Lie}}
\newcommand{\Aut}{{\rm Aut}}
\newcommand{\Out}{{\rm Out}}
\newcommand{\Inn}{{\rm Inn}}
\newcommand{\Image}{{\rm Image}}
\newcommand{\Et}{{\rm Et}}
\newcommand{\Gr}{{\rm Gr}}
\newcommand{\gr}{{\rm gr}}
\newcommand{\gl}{{\rm gl}}

\newcommand{\Id}{{\rm Id}}
\newcommand{\rk}{{\rm rk}}
\newcommand{\pardeg}{{\rm par.deg}}
\newcommand{\Res}{{\rm Res}}
\newcommand{\Fr}{{\rm Frob_p}}
\newcommand{\Spec}{{\rm Spec}}
\newcommand{\Ext}{{\rm Ext}}
\newcommand{\Sym}{{\rm Sym}}
\newcommand{\Tor}{{\rm Tor}}
\newcommand{\Supp}{{\rm Supp}}
\newcommand{\depth}{{\rm depth}}
\newcommand{\coh}{{\rm Coh}}
\newcommand{\nr}{{\rm Norm}}
\newcommand{\codim}{{\rm codim}}
\newcommand{\Gal}{{\rm Gal}}
\newcommand{\Ram}{{\rm Ram}}
\newcommand{\Mor}{{\rm Mor}}
\newcommand{\rad}{{\rm rad}}
\newcommand{\ch}{{\rm char}}
\newcommand{\cov}[1]{{\rm Cov}^{\rm top}_{#1}}
\newcommand{\ft}{ \pi^{\rm top}_1}
\newcommand{\fa}{ \pi^{\rm alg}_1}

\newcommand{\qtq}[1]{\quad\mbox{#1}\quad}

%%%%%%%%%%%%%%%%%%%%%%%%%%%%%%%%%%%%%%%%%%%%%%%%%%%%%%%%%%%%%%

\title[Calabi-Yau Varieties with Semi-stable Fibre Structures]{Calabi-Yau Varieties with Semi-stable Fibre Structures}

\author[Yi Zhang ]{Yi Zhang}
\address{\rm School of Mathematical Sciences, Fudan
University, Shanghai, 200433, People's Republic of
China}\email{zhangyi\_math@fudan.edu.cn}

\author[Kang Zuo]{Kang Zuo}
\address{\rm Faculty of Mathematics \& Computer Science,
Johannes Gutenberg Universitaet Mainz, Staudingerweg 9, D-55128
Mainz, Germany}
\email{kzuo@mathematik.uni-mainz.de}

\thanks{The first author was supported by the National
Natural Science Foundation of China(\#10731030). This paper was
supported by the joint Chinese-German project ``Komplexe
Geometrie''(DFG-NSFC)}

\maketitle \pagestyle{myheadings} \markboth{\hfill Y. Zhang and K.
Zuo \hfill}{\hfill  \hfill} \setcounter{tocdepth}{1}
\setcounter{page}{1}
\begin{abstract}
Motivated by the Strominger-Yau-Zaslow conjecture, we study
Calabi-Yau varieties with semi-stable fibre structures.  We use
Hodge theory to study the higher direct images of wedge products of
relative cotangent sheaves of  certain semi-stable families over
higher dimensional quasi-projective bases, and obtain some results
on positivity. We then apply these results to study nonisotrivial
Calabi-Yau varieties fibred by semi-stable Abelian varieties (or
hyperk\"ahler varieties).
\end{abstract}

The mirror symmetry conjecture is based on the suggestion that two
sigma-models in superstring theory are equivalent. The targets of
the models are Calabi-Yau 3-folds, so mirror symmetry predicts that
these should come in pairs, $M$ and $\breve{M}$, satisfying
$h^{p,q}(M)=h^{p,3-q}(\breve{M}).$ Recently Strominger-Yau-Zaslow
noticed in \cite{SYZ} that String Theory suggests that if $M$ and
$\breve{M}$ are mirror pairs of n-dimensional CY manifolds, then on
M should exist a special Lagrangian $n$-tori fibration $f:M\>>> B,$
(with some singular fibres) such that $\breve{M}$ is obtained by
finding some suitable compactification of the dual fibration.\\

Motivated by the Strominger-Yau-Zaslow conjecture, we study the
fibration $f: X\to Y,$ where $X$ is a Calabi-Yau type projective
manifold.
\begin{definition}(Calabi-Yau manifolds).
\begin{myenumi}

\item A projective complex manifold $Y$ is called  \emph{Calabi-Yau} if
    the following conditions are satisfied :
\begin{myenumii}
\item the canonical line bundle $\omega_Y$ of $Y$ is trivial;

\item $H^0(Y,\Omega_Y^p)=0$ for $p$ with $0<p<\dim Y.$
\end{myenumii}
\item  A compact K\"ahler manifold $Y$ is called
\emph{hyperk\"ahler} if it is
    of dimension $2n\geq 4$ and the following conditions are satisfied :
\begin{myenumii}
\item  there is a non-zero holomorphic two
    form $\beta_Y$ unique up to scalar such that $\det(\beta_Y)$ is
    nowhere zero;
    \item $H^1(Y,\sO_Y)=0.$
\end{myenumii}

\item A compact K\"ahler manifold $Y$ is called \emph{Calabi-Yau
type} if its canonical line bundle $\omega_Y$ is trivial. (Thus,
Ableian varieties and hyperk\"ahler manifolds are Calabi-Yau-type.)
\end{myenumi}
\end{definition}
\begin{myrem}If $Y$ is hyperk\"ahler, then $\dim
H^2(Y,\sO_Y)=1,$ and $\omega_X$ is trivial.
\end{myrem}

Studying projective manifolds with semi-stable fiber structures by
Deligne's Hodge theory and the semi-positivity package, we then have
:

\begin{theorem}\label{main-theorem-1}
Let $f:X\to Y$ be a semi-stable family between two non-singular
projective varieties with
$$f:X_0=f^{-1}(Y_0) \to Y_0$$ smooth, $S=Y \setminus Y_0$ a reduced normal crossing divisor and
$\Delta=f^*S$ a relative reduced normal crossing divisor in $X.$
Assume that $f$ satisfies that
\begin{myenumii}
    \item the polarized VHS $R^kf_*\Q_{X_0}$ is strictly of weight $k;$
    \item $H^{k}(X,\sO_X)=0;$
    \item $Y$ is simply connected.
\end{myenumii}
Then  $f_*\Omega^k_{X/Y}(\log \Delta)$ is locally free on $Y$
without flat quotient, and $S\neq \emptyset.$ Moreover,
$$\deg_C (f_*\Omega^k_{X/Y}(\log \Delta))>0$$ for any {\it
sufficiently general} curve $C\subset \overline{M}.$
\end{theorem}

Based on the above theorem, we consider Calabi-Yau varieties fibred
by semi-stable Abelian varieties (or by hyperk\"ahler varieties)
over $\C\P^1,$ and then we obtain the following one of our main
results:

\begin{theorem}\label{main-theorem-3}Let $f:X\to \P^1$ be a semi-stable family fibred by Abelian varieties such that
$$f:X_0=f^{-1}(C_0) \to C_0$$
is smooth with finite singular values $S=\P^1 \setminus C_0$ and a
normal crossing $\Delta=f^{-1}(S).$ Assume that $X$ is a projective
manifold with  trivial canonical sheaf $\omega_X$ and
$H^1(X,\sO_X)=0.$ Then, $f$ is nonisotrivial and $\dim X \leq 3.$

In particular, if $X$ is a Calabi-Yau manifold then $X$ is one of
the following:
\begin{myenumii}
  \item $K3$ with $\#S \geq 6$. %By well known results of Kodaira, $\#S = 24$.
  Moreover, if $\#S =6$ then $X\to \P^1$ is modular, i.e., $C_0$ is the quotient of
the upper half plane $\sH$ by a subgroup of $\mathrm{SL}_2(\Z)$ of
finite index.
    \item Calabi-Yau threefold with $\#S \geq 4$. Moreover, if $\# S =4$ then this family is rigid and there exists an \'etale covering $\pi:Y'\to
\P^1$ such that $f':X'=X\times_{\P^1} Y'\to Y'$ is isogenous over
$Y'$ to a product $ E\times_{Y'} E ,$ where $h: E\to Y'$ is a family
of semi-stable elliptic curves and
modular.%i.e., $C_0$ is the quotient of the upper half plane $\sH$
%by a subgroup of $\mathrm{SL}_2(\Z)$ of finite index.

\end{myenumii}
\end{theorem}

Furthermore, we consider Calabi-Yau varieties fibred by semi-stable
Abelian varieties with higher dimension base $Y,$ the result in
\cite{ZQ} shows that the base $Y$ must be rational connected(cf.
Fact \ref{observation-2} in section 2).

Consider a rationally connected manifold $Z.$ Let $D_\infty$ be a
reduced normal crossing divisor in $Z.$  One can choose a very
freely rational curve (not necessary smooth) $C$ intersecting each
component of $D_\infty$ transversely, then
$$\pi_1(C\cap (Z-D_\infty)) \onto \pi_1(Z-D_\infty)\to 0 \mbox{ is surjective}.$$
Actually, Koll\'ar's results on fundamental groups show the
following fact:

\begin{proposition}[cf. \cite{Ko2000}]\label{Kollar-fundamental-grp}%[A quasi-projective version of the Lefschetz hyperplane theorem]
Let $X$ be a smooth projective variety and $U\subset Z$ be an open
dense subset such that $Z\setminus U$ is a normal crossing divisor.
Assume that $Z$ is rationally connected. Then, there exists a very
free rational curve $C\subset Z$ such that it intersects each
irreducible component of $Z\setminus U$ transversally and the map of
topological fundamental groups $ \pi_1(C\cap U) \onto \pi_1(U)\to 0$
is surjective; moreover

\begin{myenumii}

\item  if $\dim Z \geq 3$ then $C$ is a smooth rational curve in
$Z,$ and

\item  if $\dim Z =2$ then there is an immersion $h: \P^1 \to C \subset Z.$

\end{myenumii}
\end{proposition}

With the above proposition \ref{Kollar-fundamental-grp}, we apply
Theorem \ref{main-theorem-3} and immediately obtain the following
result:

\begin{theorem}\label{main-theorem-4}
Let $f:X\to Y$ be a semi-stable family of Abelian varieties between
two non-singular projective varieties with $f:X_0=f^{-1}(Y_0) \to
Y_0$ smooth, a reduced normal crossing divisor $S=Y \setminus Y_0$
and a relative reduced normal crossing divisor $\Delta=f^{*}(S)$ in
$X.$ Assume that
\begin{myenumii}
\item the period map of the VHS $R^1f_*(\Q_{X_0})$
is injective at one point in $Y_0;$  %ing induced moduli morphism is generally finite and Assume that
\item the canonical bundle $\omega_X$ is trivial and
$H^0(X,\Omega^1_X)=0.$
\end{myenumii}
Then, the dimension of a general fibre is bounded above by a
constant  dependent on $Y.$
\end{theorem}

With similar methods, we consider any Calabi-Yau variety fibred by
semi-stable hyperk\"ahler varieties, and show the base manifold is
rationally connected such that the dimension of a general fibre is
bounded above by a constant depending on the base, moreover if the
base is $\C\P^1$ then the dimension of a general fibre is four (see
Theorem \ref{dim-hyperkahler} and Theorem \ref{main-theorem-5} in
section 3).

%\newpage

{\small\tableofcontents}

\section{Preliminary}
%\section{Notations and basic results}

Let $M$ be a smooth $m$-dimensional complex quasi-projective
variety, $M\subset\overline{M}$ a smooth projective compactification
such that $\overline{M}-M=D_\infty$ is a reduced normal crossing
divisor. Let $j:M\>\subset>> \overline{M}$ be the open embedding.

Let $\mathbf{H}$ be an arbitrary polarized variation of $\R$-Hodge
structures (VHS for short) over $M$ with unipotent monodromies
around $D.$ Let $\sH=\mathbf{H}\otimes \sO_M,$ a holomorphic vector
bundle, $\sH=\sF^0\supset\sF^1\supset\cdots\supset\sF^w\supset 0$
the Hodge filtration with vector bundles(Hodge bundles). $w$ will be
called the weight of the VHS. $\sF^w$ will be denoted $\sF^b,$ here
the sign $b$ is for bottom, i.e., $\sF^b$ is the lowest piece of the
Hodge filtration(cf. \cite{Sch}).

\subsection*{Deligne's canonical extension of a VHS}
We consider a \textbf{special coordinate} around $o \in D,$ i.e., it
is a coordinate neighborhood $\triangle_o\subset \overline{M}$ of
$o$ which is isomorphic to $\Delta^n,$ and
$$M\cap \triangle_o=\{z=(z_1,\cdots, z_n)\,\,|\,\, z_1\neq 0,\cdots, z_l\neq 0\mbox{ for some } 1\leq l\leq n\}\cong (\Delta^*)^l\times \Delta^{n-l}.$$
In particular, $$\triangle_o\cap D_\infty=\{z=(z_1,\cdots,
z_n)\,\,|\,\, z_1\cdots z_l=0\}.$$

Let $\gamma_\alpha$ be a local monodromy around $z_\alpha=0$ in
$\triangle_o$ for $\alpha=1, \cdots l.$ $N_\alpha=\log
\gamma_\alpha$ is thus nilpotent.

Choose a flat multivalued basis $(v_{\cdot})$ of $\sH$ over
$\triangle_o\cap M.$  The formula
$$(\widetilde{v_{\cdot}})(z):=\exp(\frac{-1}{2\pi\sqrt{-1}}\sum_{\alpha=1}^l\log z_\alpha N_\alpha )(v_{\cdot})(z)$$
gives a single-valued basis basis of $\sH.$ Deligne's canonical
extension $\overline{\sH}$ of $\sH$ to $\triangle_o$ is generated
by $(\widetilde{v_\cdot})$ (cf. \cite{Del}, \cite{Sch}). The
construction of $\overline{\sH}$ is independent of the choice of
$z_i's$ and $(v_\cdot).$ Obviously, $\overline{\sH}$ is locally
free.

Denote by $\overline{\sF}^p:=\overline{\sH}\cap j_*\sF^p.$ Thanks to
the nilpotent orbit theorem
(cf.\cite{CKS},\cite{Sch}),$$\overline{\sH}=\overline{\sF}^0\supset\overline{\sF}^1\supset\cdots\supset\overline{\sF}^w\supset
0$$ is a  filtration of locally free sheaves.

\subsection*{The Higgs bundle associated to a VHS}
Let $\nabla$ be the Gauss-Manin connection on the Hodge bundle
$\sH.$ The Deligne canonical extension $\overline{\sH}$ of $\sH$ and
the nilpotent theorem allow the  Gauss-Manin connection extend to be
an regular connection
$$\overline{\nabla}: \overline{\sH} \rightarrow \overline{\sH}
\otimes \Omega^1_{\overline{M}}(\log D_\infty)$$ with the Griffiths
transversality:
$$\overline{\nabla}: \overline{\sF}^p \rightarrow
\overline{\sF}^{p-1} \otimes \Omega^1_{\overline{M}}(\log
D_\infty)\,\, p=1,\cdots, w.
$$

Denote by $$\overline{E}^{p,w-p}=\overline{\sF}^p/
\overline{\sF}^{p+1} \,\, p=1,\cdots, w.$$ the Griffiths
transversality allows $\sO_{\overline{M}}$-linear morphisms
$$\theta^{p,q}: \overline{E}^{p,q}\rightarrow
\overline{E}^{p-1,q+1}\otimes \Omega^1_{\overline{M}}(\log
D_\infty),$$ for $p=1,\cdots, w$ and $p+q=w.$ We then have the
Higgs bundle
$$(\overline{E}:=\bigoplus \overline{E}^{p,q}, \theta:=\bigoplus \theta_{p,q}),$$
i.e., $\theta :\overline{E}\rightarrow \overline{E}\otimes
\Omega^1_{\overline{M}}(\log D_\infty)$ is an
$\sO_{\overline{M}}$-liner morphism with $\theta\wedge\theta=0$ (cf.
\cite{Si1}).

\subsection*{Koll\'ar's decomposition}
We recite the results of Proposition 4.10 and Remark 4.11 in
\cite{Kol}.

\begin{definition}\label{sufficiently-general-curve}
A smooth projective curve $C \subset \overline{M}$ is {\it
sufficiently general} if it satisfies that
\begin{enumerate}
    \item $C$ intersects
$D_\infty$ transversely;%(then $S=C\cap D_{\infty}$ is a finite set);
    \item $\pi_1(C_0) \onto \pi_1(M)\to 0 \mbox{ is surjectve}$ where $C_0=C\cap
M.$
\end{enumerate}
\end{definition}
\begin{myrem}There are many sufficiently general curves :
Let $C$ be a complete intersection of very ample divisors such that
it is a smooth projective curve
in $\overline{M}$ intersecting $D_\infty$ transversally. %The homomorphism
The quasi-projective version of the {\it Lefschetz hyperplane
theorem}(cf.\cite{GM}) guarantees the subjectivity of
$\pi_1(C_0)\onto \pi_1(M).$ But we note that a curve $C$ is
sufficiently general does not mean that it must be a complete
intersection of hyperplanes.
\end{myrem}

\begin{theorem}[Kollar]\label{Kollar-decomposition}
Let $M$ be a quasi-projective $n$-fold with  a smooth projective
completion $\overline{M}$ such that $D_{\infty}=\overline{M}-M$ is a
reduced normal crossing divisor. Consider the lowest piece $\sF^b$
of the Hodge filtration of a polarized VHS $H$ on $M.$ Assume that
all local monodromies of $\V$ are unipotent. There is a unique
decomposition for the canonical extension $\overline{\sF}^b$ of
$\sF^b$
$$\overline{\sF}^b=\sA \oplus \sU,$$ such that
$\sA$ has no flat quotient and  $\sU$ is a unitary flat bundle on
$\overline{M}.$ Moreover,if $C$ is a {\it sufficiently general}
curve in $\overline{M}$ then $\sA|_C$ is an ample vector bundle on
$C.$
\end{theorem}

\subsection*{The Fujita-Kawamata's Semi-positive package}
We recite the Semi-positive package in \cite{Kaw} and \cite{Kaw1}.
\begin{definition}\label{URC}
Let $\pi:X\to Y$ be an algebraic fibration  with $d=\dim X-\dim Y.$
We say $\pi$ satisfies the {\it unipotent reduction condition} ({\it
URC}) if the following are held :
\begin{enumerate}
\item there is a Zariski open dense subset $Y_0$ of $Y$ such that
$D=Y\setminus Y_0$ is a {\it divisor of normal crossing on $Y$},
i.e., $D$ is a reduced effective divisor and if $D=\sum_{i=1}^N
D_i$ is the decomposition to irreducible components, then all
$D_i$ are non-singular and cross normally;

\item $\pi: X_0 \to Y_0$ is smooth where $X_0=\pi^{-1}(Y_0);$

\item all local monodromies of $R^d\pi_{*}\Q_{X_0}$ around $D$ are
unipotent.
\end{enumerate}
\end{definition}
The ${\it URC}$ holds automatically for any semi-stable family.

\begin{theorem}[Kawamata]\label{kawamata}
Let $\pi:X\to Y$ be a proper algebraic family with connected fibre
and $\omega_{X/Y}:=\omega_{X}\otimes \pi^*\omega^{-1}_Y$ be the
relative dualizing sheaf.  Let $\sF$ be the bottom filtration of the
VHS $R^d\pi_{*}\Q_{X_0}$ where $X_0=f^{-1}(Y_0)$ and $d=\dim X-\dim
Y.$ Assume that $\pi$ satisfies {\it URC}, then we have :
\begin{myenumi}
    \item $\pi_*\omega_{X/Y}=\overline{\sF}$ is locally
 free,where $\overline{\sF}$ is the Deligne canonical extension of $\sF$ over $\overline{M}.$
    \item $\pi_*\omega_{X/Y}$ is {\it semi-positive}, i.e., for every projective curve $T$ and
morphism $g: T \to Y$ every quotient line bundle of
$g^*(\pi_*\omega_{X/Y})$ has non-negative degree.
\end{myenumi}
\end{theorem}

In this paper, we also need the following  result in  \cite{Kaw2}:
\begin{corollary}[Kawamata]\label{big-torelli} Let $Y$ be projective
manifold and let $Y_0$ be a dense open set of $Y$ such that
$S=Y\setminus Y_0$ is a reduced normal crossing divisor. Let
$\mathbf{H}$ be a polarized VHS of strict weight $w$ over $Y_0$ such
that all local monodromies are unipotent.

Assume that
$$T_{Y,p} \> >> \mathrm{Hom}(\sF_{p}^w, \sF^{w-1}_p/\sF^w_p) $$
is injective at $p\in Y_0$ where $\sF^w=: \sF^b$ is the lowest
piece of the Hodge filtration of the VHS $\mathbf{H}.$ Then, $\det
\overline{\sF}^b$ is a {\it big} line bundle,where
$\overline{\sF}^b$ is the canonical extension of $\sF^b.$
\end{corollary}
\begin{proof}%[Sketch of the proof of \ref{big-torelli}]
The form $c_1(\sF^b,h)$ on $Y_0$ is a current on $Y$ and represents
the Chern class $C_1(\overline{\sF}^b)$ due to Cattani-
Kaplan-Schmid's theorem (cf.\cite{CKS}).We can calculate
$c_1(\sF^b,h)$ as follows :

Let $\sN=(\sF^b)^{\vee}$ and $h$ is the Hodge metric on
$\sH=\mathbf{H}\otimes \sO_{Y_0}.$ Then,  we have the curvature form
$$ \Theta(\sN,h)=-\theta\wedge\bar\theta_h|_{\sN}+\bar A_h\wedge A
\mbox{ on } Y_0
$$ where $A$ is the second fundamental form of $\sN$ (cf. \cite{Sch}), and $\theta(\sN)=0$ on $Y_0$ by the
Griffiths transversality.

On the other hand, the lemma 2.6 in \cite{V83} or Theorem 0.1 in
\cite{Zuo0} shows that  $$\int_{Y_0} c_1(\sF^b,h)^{\dim Y} > 0,$$
by the injectivity of the morphism $T_{Y,p}
\to\mathrm{Hom}(\sF_{p}^w, \sF^{w-1,1}_p).$ Since
$\det(\overline{\sF}^b)$ is  a {\it nef} line bundle by Theorem
\ref{kawamata}, The Rienamnn-Roch Theorem says that
$\det(\overline{\sF})$ is {\it big}, or we can obtain the bigness
by Sommese-Kawamata-Siu's numerical criterion:

\emph{If $L$ is a hermitian semi-positive line bundle on a compact
complex manifold $X$ such that $\int_X \wedge^{\dim X}c_1(L)>0,$
then $L$ is a {\it big} line bundle on $X.$ In particular, if the
line bundle $L$ is {\it nef } with $(L)^{\dim X}>0,$ then $L$ is
{\it big}.}\\
\end{proof}

\subsection*{Families of Calabi-Yau-type manifolds}

Let $(Y,L)$ be a polarized manifold. Let
$$\pi:(\ssY,(Y,L)) \to (\mathfrak{N}_{c_1(L)},0)$$ be the
polarized Kuranishi family and $D$  the classifying space of the
polarized $\Q$-Hodge structure of
$$H_{prim}^n(Y,\C):=\Ker(H^n(Y,\C) \> \wedge c_1(L) >>
H^{n+2}(Y,\C)).$$
\begin{theorem}[Bogomolov-Tian-Todorov \cite{Ti},\cite{To89}]
Assume that $Y$ has trivial canonical line bundle $\omega_Y.$
Then, the Kuranishi family of $(Y,L)$ is universal and the base
manifold is a smooth open set in the Euclidean space of dimension
$$\dim_\C H^1(Y,\Theta_Y)_{c_1(L)},$$ where
$H^1(Y,\Theta)_{c_1(L)}:=\Ker(H^1(Y,\Theta_Y) \> \wedge c_1(L) >>
H^2(Y,\sO_Y)).$
\end{theorem}

Let  $$f:(\sY,(Y,L)) \to (Z,0)$$ be a family of polarized manifolds.
Locally, near $0\in Z,$ one has a commutative diagram of period maps
$$
\begin{TriCDV}
{\mathfrak{N}_{c_1(L)}} {\< \iota <<} {Z} {\SE \Phi EE}{\SW W \Phi_Z
W} {D}
\end{TriCDV}
$$

Let $\mu_0=(d\Phi)_0$ and $\lambda_0=(d \Phi_{Z})_0.$ Also we have
the Kodaira-Spencer map $\rho=(d\iota)_0.$

\begin{theorem}[Griffiths's Infinitesimal Torelli Theorem, \cite{Gr70}]
We have a natural morphism
$$\mu_0 : T_{\mathfrak{N}_{c_1(L),0}}=H^1(Y,\Theta_Y)_{c_1(L)}\>>> \bigoplus_p \Hom(H^{n-p,p}_{prim},H^{n-p-1,p+1}_{prim})$$
and $\lambda_0=\mu_0 \circ \rho,$ where
$T_{\mathfrak{N}_{c_1(L),0}}$ is the Zariski tangent space. Write
$\mu_0=(\mu_0^0,\cdots,\mu_0^n)$ where
$$\mu^i_0 : T_{\mathfrak{N}_{c_1(L),0}}=H^1(Y,\Theta_Y)_{c_1(L)}\>>> \Hom(H^{n-i,i}_{prim},H^{n-i-1,i+1}_{prim}).$$
If $Y$ has trivial $\omega_Y$ then $\mu^0_0$ is an isomorphism and
$\mu_0$ is an injective map.

\end{theorem}

In this paper, we will use the following fact frequently: \emph{For
a non-isotrivial family of manifolds with trivial canonical line
bundle, the period map is not degenerate at a general point of the
base.} The Infinitesimal Torelli Theorem shows that \emph{the
condition that the period map $\Phi$ is not degenerate at $0$ is
equivalent to the condition that the Kodaira-Spencer map
is injective at $0.$}\\

\section{Calabi-Yau type manifolds with semi-stable fibrations.}

From now on, we always assume that $X$ is a K\"ahler manifold with
trivial canonical line bundle $\omega_X\cong \sO_X.$\\

\subsection*{Some Facts}

Let $f:X \to C$ be a semi-stable family from a projective manifold
$X$ with trivial $\omega_X$ to a smooth projective curve $C.$ Since
$$\omega_C^{-1}=\sO_C(\sum t_i-\sum t_j)\mbox{ with }
\#\{i\}-\#\{j\}=2-2g(C),$$ we then have
$$\omega_{X/C}=\omega_{X}\otimes
f^*\omega_C^{-1}=f^*\omega_C^{-1}=\sO_{X}(\sum X_{t_i}-\sum
X_{t_j}).$$ Thus, $f_*\omega_{X/C}=\omega_C^{-1}$ and
$\sO_{X}(X_{t_i})|_{X_t}=\sO_{X_t}(X_{t_i}\cdot X_t)=\sO_{X_t} \,
\forall t\in C$ where $X_t=f^{-1}(t).$ Hence,
$$\omega_{X_t}=\omega_{X/C}|_{X_t}=\sO_{X_t} \forall t\in C.$$
i.e.,  the canonical sheaf of each closed fibre is trivial.

On the other hand, if all closed fibres have trivial canonical
bundle then $\omega_X$ is trivial if and only if
$f_*\omega_{X/C}=\omega_C^{-1}.$ In particular, if the total space
is a Calabi-Yau manifold, then the fibres can be Abelian varieties,
lower dimensional Calabi-Yau varieties or hyperk\"ahler varieties.

\begin{fact}\label{big-nef-CY} Let $f: X\to Y$ be a semi-stable  family of
Calabi-Yau varieties over a higher dimension base such that $f$ is
smooth over $Y_0$ and $Y\setminus Y_0$ is a reduced normal
crossing divisor. We have :
\begin{myenumi}
\item If the induced moduli map is a generically finite
morphism, then the line bundle $f_*\omega_{X/Y}$ is {\it big and
nef}.

\item If $f$ is a smooth family and the induced period
map has no degenerated point, then the line bundle $f_*\omega_{X/Y}$
is ample.
\end{myenumi}
\end{fact}

%\begin{proof}Let $X_0=f^{-1}(Y_0)$ (resp. $X$) and $h$ be the induced Hodge metric on
%$f_*\omega_{X_0/Y_0}$ (resp. $f_*\omega_{X/Y}$). The
%Bogomolov-Tian-Todorov theorem shows that
%$$c_1(f_*\omega_{X_0/Y_0},h)=\omega_{WP}$$ where $\omega_{WP}$ is
%the K\"ahler form of the Weil-Petersson metric over $Y_0.$
%\begin{myenumi}

%\item $\omega_{WP}$ is a positive $(1,1)$ form on a Zariski open
%subset of $Y_0.$ The Cattani-Kaplan-Schmid theorem shows that the
%Chern form $c_1(f_*\omega_{X_0/Y_0},h)$ on $Y_0$ is a current and
%represent the Chern class $c_1(f_*\omega_{X/Y})$ on $Y.$
%$f_*\omega_{X/Y}$ is {\it big} line bundle because it is {\it nef}
%and
%$$\int_{Y_0} c_1(f_*\omega_{X_0/Y_0},h)^{\dim Y} >0 .$$

%\item $c_1(f_*\omega_{X/Y},h)=\omega_{WP} $ is a positive $(1,1)$
%form on $Y$ because the period map has no degenerated point. The
%second statement follows from the Kodaira embedding theorem
%directly : {\it A compact complex manifold $M$ is an algebraic
%variety if and if only it has a closed, positive $(1,1)$-form
%$\omega$ whose cohomology class $[ \omega ]$ is rational.}
%\end{myenumi}
%\end{proof}
In particular, we have:
\begin{corollary}\label{observation-1}Let $f: X\to C$ be a semi-stable
non-isotrivial family over a smooth projective curve $C,$ where $X$
has trivial $\omega_X.$ Then, $f_*\omega_{X/C}$ is a big line bundle
and $C$ is $\P^1.$
\end{corollary}

%We can also deduce this observation by  the following proposition
%\ref{corollary-Kollar-decomposition}.\\

Let $f:X \to Y$  be a semi-stable proper family smooth over a
Zariski open dense set $Y_0$ such that $S=Y-Y_0$ is a reduced normal
crossing divisor. Suppose that $X$ is a projective manifold. Then, a
general fibre is a smooth projective manifold  with trivial
canonical bundle and it has certain geometric type `K'. We know that
the coarse quasi-projective moduli scheme $\mathfrak{M}_K$ exists
for the set of all polarized projective manifolds with trivial
canonical sheaf and with given type `K' (cf.\cite{V95}).

By the infinitesimal Torelli theorem, that the family $f$
satisfies the condition in the corollary \ref{big-torelli} if and
only if the unique moduli morphism  $\eta_f: Y_0 \to
\mathfrak{M}_K$ for $f$ is a generically finite morphism.
Moreover, the condition is equivalent to that $f$ contains no
isotrivial subfamily whose base is a subvariety passing through a
general point of $Y.$ If $Y$ is a curve, that $f$ satisfies the
condition in the corollary \ref{big-torelli} if and only if the
family $f$ is non-isotrivial. Actually, we have a result from
Koll\'ar's decomposition and the semi-positive package:

\begin{corollary}\label{corollary-main-theorem-0}
Let $f:X\to Y$ be a surjective morphism between two non-singular
projective varieties such that every fibre is irreducible and $f:
X_0=X\setminus \Delta \to Y\setminus S$ be the maximal smooth
subfamily where $S=Y \setminus Y_0$ is a reduced normal crossing
divisor. Let $\sF$ be the lowest piece $\sF^b$ of the Hodge
filtration of a polarized VHS $R^{n-1}f_*(\Q_{X_0})$ on $Y_0,$ i.e.,
$$\sF=F^{n-1}(R^{n-1}f_*(\Q_{X_0})\otimes \sO_{Y_0}).$$

Assume that $X$ is a projective $n$-fold with trivial canonical
sheaf $\omega_X$  and $f$ is {\it semi-stable}, i.e., $\Delta=f^*S$
is a relative  reduced normal crossing divisor in $X.$.

If the moduli morphism of $f$ is generically finite, then we have :

\begin{myenumi}
\item $f_*\Omega^{n-1}_{X/Y}(\log \Delta)=f_*\omega_{X/Y}=\overline{\sF}$
is  a {\it big} and {\it nef} line bundle on $Y,$ where
$\overline{\sF}$ is  the canonical extension of $\sF.$

\item Moreover, for any {\it sufficient general} curve $C\subset
Y$  $\overline{\sF}|_C$ is ample on $C.$
\end{myenumi}
\end{corollary}

Recently, %using the theorem of Graber-Harris-Starr (cf.\cite{GHS}) and Viehweg's {\it weakly positivity},
Zhang(cf.\cite{ZQ}) proves that $\log Q$-Fano varieties are
rationally connected, and it implies a result which was obtained by
Koll\'ar-Miyaoka-Mori in case of threefold (cf.\cite{KMM92b}): {\it
Any higher dimensional variety with a {\it big} and {\it nef}
anticanonical bundle must be rationally connected}.

Therefore, we have:
\begin{fact}\label{observation-2}
Let $f: X\to Y$ be a semi-stable proper family between two
non-singular projective varieties. Assume that $X$ has trivial
canonical sheaf $\omega_X$ and the induced moduli morphism $\eta_f$
is generically finite. Then, the anti-canonical line bundle
$\omega_Y^{-1}$ of $Y$ is {\it big} and {\it nef}, and so $Y$ is
{\it rationally connected}.
\end{fact}

\subsection*{Vanishing of unitary flat subbundles}
Now we consider the cohomology geometry in case of Calabi-Yau
manifolds. Before proving Theorem \ref{main-theorem-1}, we recite
Deligne's Theory as follows:

Let $f:X\to Y$ be a semi-stable family between two non-singular
projective varieties with
$$f:X_0=f^{-1}(Y_0) \to Y_0$$ smooth, $S=Y \setminus Y_0$ a reduced normal crossing divisor and
$\Delta=f^*S$ a relative reduced normal crossing divisor in $X.$ Let
$y\in Y_0$ be a fixed point. Consider the following commutative
diagram:
$$
\begin{TriCDV}
{H^m(X,\Q)}{\> i^*
>>}{H^m(X_0,\Q)} {\SE  \overline{i}_{y}^* EE}{\SW  W i_{y}^* W}
{H^m(X_{y},\Q)}
\end{TriCDV}
$$
where $i_{y}:X_y\hookrightarrow X_0,$
$\overline{i}_{y}:X_y\hookrightarrow X$ are natural embeddings.
$H^m(X_0,\Q)$ can be equipped with a functorial mixed Hodge
structure $(W_m,F^p,H^m(X_0,\Q))$ of pure weight $m$ with
$$W_m(H^m(X_0,\Q))=\Image(H^m(X,\Q) \> i^* >> H^m(X_0,\Q)).$$
 In
addition, $H^m(X,\Q)$ $H^m(X_y,\Q)$ also have a pure Hodge structure
structure of weight $m, $ and all of $i^*, i_y^*,
\overline{i}_{y}^*$ are morphisms of mixed Hodge structures. It
follows that $$\Image (\overline{i}_{y}^*) = \Image (i_{y}^*).$$(cf.
Section 3, in particular Corollary 3.2.18 of \cite{Del}).

$$i_y^*:H^m(X,\Q)\>>>H^m(X_y,\Q)$$
is a morphism of mixed Hodge  structures.

We note that $\pi_1(Y_0,y)$ acts on the $H^n(X_y,\Q).$ In
particular, there is a $\pi_1$-invariant subspace
$$H^m(X_y,\Q)^{\pi_1}\hookrightarrow H^m(X_y,\Q).$$ We glue
these invariant subspaces together into a constant sheaf $\sI$ on
$Y_0.$ $\sI$ coincides with the constant sheaf on $Y_0$ with the
fiber $H^0(Y_0,R^mf_*\C_{X_0}),$ which is a subsheaf of
$R^mf_*\Q_{X_0}$ generated by the global sections of
$R^nf_*\Q_{X_0},$ i.e.,
$$\sI=(R^mf_*\Q)^{\pi_1}.$$ There is
a natural isomorphism of $\Q$-vector spaces
$$\kappa_y: H^0(Y_0,R^mf_*\Q_{X_0}) \> \simeq >> H^m(X_y,\Q)^{\pi_1}.$$
Deligne shows that Leray's spectral sequence for $f: X_0\to Y_0$
$$E^{p,q}_2=H^p(Y_0,R^qf_*\Q_{X_0})\Rightarrow H^{p+q}(X_0,\Q)$$
degenerates at $E_2$(cf. Theorem 4.1.1 (i) of \cite{Del}), and so it
follows that the canonical mapping
$$j:H^m(X_0,\Q)\>>> H^0(Y_0,R^mf_*\Q_{X_0})\>>>0$$ is
surjective. Moreover,  $i_y^*$ and $\overline{i}_{y}^*$ can be
decomposed as follows :
\begin{equation*}
i^*_y: H^m(X_0,\Q)\> j >>  H^0(Y_0,R^mf_*\Q_{X_0}) \> \kappa_y >>
H^m(X_y,\Q)^{\pi_1}\hookrightarrow H^m(X_y,\Q),
\end{equation*}
\begin{equation*}
\overline{i}_{y}^*:  H^m(X,\Q)\>i^*>>H^m(X,\Q)\> j >>
H^0(Y_0,R^mf_*\Q_{X_0}) \> \kappa_y >>
H^m(X_y,\Q)^{\pi_1}\hookrightarrow H^m(X_y,\Q).
\end{equation*}
Since $i_y^*, \overline{i}_{y}^*$  are morphisms of mixed Hodge
structure and $$\Image (\overline{i}_{y}^*) = \Image
(i_{y}^*)=H^m(X_y,\Q)^{\pi_1},$$ $H^m(X_y,\Q)^{\pi_1},$ as an image
of a Hodge structure, is a Hodge substructure in $H^m(X_y,\Q).$ (cf.
the proof of Theorem 4.11 (ii) of \cite{Del}).

We can introduce a Hodge structure on $H^0(Y_0,R^mf_*\Q_{X_0}),$ as
a quotient structure on $H^m(X_0,\Q).$ This quotient structure is
independent of $y\in Y_0$(cf. Corollary 4.1.2 of \cite{Del}).

\begin{proof}[Proof of Theorem \ref{main-theorem-1}]
At first, we note that $f_*\Omega^k_{X/Y}(\log \Delta)$ is just the
canonical extension of $F^{k}(R^{k}f_*(\Q_{X_0})\otimes \sO_{Y_0})$
which is the lowest piece of the Hodge filtration of the VHS
$R^kf_*\Q_{X_0},$ thus $f_*\Omega^k_{X/Y}(\log \Delta)$ is locally
free(cf. Lemma 1 in \cite{Kaw1} for detail). That the polarized VHS
$R^kf_*\Q_{X_0}$ is of weight $k$ guarantees that
$f_*\Omega^k_{X/Y}(\log \Delta) \neq 0.$

By Kollar's decomposition \ref{Kollar-decomposition}, we have :
$$f_*\Omega^k_{X/Y}(\log \Delta)=\sA\oplus\sU$$ such that $\sA$ has no flat quotient and
$\sU$ is flat (so $\sU$ is trivial here since $Y$ is simply
connected). We show that $f_*\Omega^k_{X/Y}(\log \Delta)$ has no
flat direct summand : Otherwise, there is a nonzero global section
$s \in H^0(Y_0,R^kf_*(\C))$ of $(k,0)$-type.

Consider Deligne's Hodge theory, we have the following commutative
diagram:
$$
\begin{TriCDV}
{H^m(X,\C)}{\> i^*
>>}{H^m(X_0,\C)} {\SE  \overline{i}_{y}^* EE}{\SW  W i_{y}^* W}
{H^m(X_{y},\C)}
\end{TriCDV}
$$
where $i_{y}:X_y\hookrightarrow X_0,$
$\overline{i}_{y}:X_y\hookrightarrow X$ are natural embeddings.
 For each pair $(p,q)$ with
$p+q=m,$ the following restriction map induced by $X_y \subset X $
is a Hodge morphism:
$$
r^{p,q}_y:H^q(X,\Omega^{p}_{X})\>\hookrightarrow >> H^m(X,\C) \>
\overline{i}_y^* >>H^m(X_y,\C)\rightarrow H^q(X_y,\Omega^{p}_{X_y}),
$$
where $y\in Y_0$ is a fixed point.

Since $$\Image (\overline{i}_{y}^*) = \Image
(i_{y}^*)=H^m(X_y,\Q)^{\pi_1(Y_0,y)},$$ we then have:

{\it The component of type $(p,q)$ of the group
$H^m(X_{y},\C)^{\pi_1(Y_0,y)}$ is just the image of
$H^q(X,\Omega^p_{X})$ under $r_y^{p,q}$.}

Let $m=k.$ We then have a nonzero lifting $\widetilde{s} \in
H^0(X,\Omega_X^k)$ of $s,$ it is a contradiction. Similarly, we
have $S\neq \emptyset.$
\end{proof}

\begin{corollary}\label{application-of-kollar}
Let $f:X\to Y$ be a semi-stable family between two non-singular
projective varieties with $$f:X_0=f^{-1}(Y_0) \to Y_0$$ smooth,
$S=Y \setminus Y_0$ a reduced normal crossing divisor and
$\Delta=f^{-1}(S)$ a relative reduced normal crossing divisor in
$X.$  Assume that $X$ is a Calabi-Yau $n$-fold and $Y$ is simply
connected. Then,

\begin{myenumii}
\item  $f$ is a nonisotrivial family with $S \neq \emptyset;$

\item $f_*\omega_{X/Y}$ is an ample line bundle over any {\it
sufficiently general} curve.
\end{myenumii}
\end{corollary}
\begin{proof}
Consider the polarized VHS $R^{n-1}f_*(\Q_{X_0}).$ Suppose that
$f$ is isotrivial then the holomorphic period map for the VHS
$R^{n-1}f_*(\Q_{X_0})$ is constant over $Y_0$ by the {\it
infinitesimal Torelli theorem}. Then the line bundle
$f_*\omega_{X_0/Y_0}$ is unitary flat over $Y_0.$ Since $f$ is
semi-stable, all local monodromies around the VHS are unipotent,
and so $f_*\omega_{X/Y}$ is unitary flat. Since $Y$ is
simply-connected, $f_*\omega_{X/Y}$ then is a trivial line bundle
on $Y,$ it is a contradiction to the theorem \ref{main-theorem-1}.
By similar arguments we then obtain  $S\neq \emptyset.$
\end{proof}
\begin{myrem}
Without assumption that $X$ is Calabi-Yau manifold and $Y$ is
simply connected. When only $X$ has trivial $\omega_X,$ we still
have the following results:
\begin{myenumi}
\item By infinitesimal Torelli theorem, if $f$ is isotrivial then
the line bundle $f_*\omega_{X/Y}$ is unitary flat .

\item Conversely, if $f_*\omega_{X/Y}$ is unitary flat and the global
Torelli theorem hold for a general fiber (e.g. a  general fiber is
K3 or Abelian variety) then $f$ is isotrivial.
\end{myenumi}
\end{myrem}

\begin{corollary}\label{kollar-flat-ample}Let $f:X\to \P^1$ be a semi-stable family with
$$f:X_0=f^{-1}(C_0) \to C_0$$ smooth, $S=\P^1 \setminus C_0,$ and
$\Delta=f^{*}S$. Assume that  $X$ is a projective manifold with
$H^{k}(X,\sO_X)=0$ and the polarized VHS $R^kf_*\Q_{X_0}$ is
strictly of weight $k.$  Then, $S \neq \emptyset$ and
$f_*\Omega^k_{X/\P^1}(\log \Delta)$ is an ample bundle on $\P^1.$
\end{corollary}

\begin{proposition}\label{character-nonisotrivial}
Let $f:X\to C$ be a semi-stable family over a smooth projective
curve with $f:X_0=f^{-1}(C_0) \to C_0$ smooth, $S=C \setminus C_0,$
and $\Delta=f^{-1}(S).$ If $X$ is a projective $n$-fold with trivial
 $\omega_X,$ then the following conditions are
equivalent :
\begin{enumerate}
    \item  $f$ is a non-isotrivial family.
    \item  $f_*\omega_{X/C}$ is an ample line bundle on $C.$
    \item  $C=\P^1$ and $\#S\geqq 3.$
\end{enumerate}
\end{proposition}
\begin{proof}
At first, we note that each smooth closed fibre of $f$ has trivial
canonical sheaf, and
$$2-2g(C)=\deg\omega_C^{-1}=\deg f_*\omega_{X/C}\geq 0.$$
\begin{myenumii}

\item $(1)\Longleftrightarrow (2)$ is proven as follows:
\begin{itemize}
\item The infinitesimal Torelli theorem holds for the VHS
$R^{n-1}f_*(\Q_{X_0})$ because there is an isomorphism for any $t\in
C_0 $ :
$$ H^1(X_t,T_{X_t})\rightarrow
\Hom(H^0(X_t,\Omega^{n-1}_{X_t}),H^1(X_t,\Omega^{n-2}_{X_t})).$$

\item $f  \mbox{ is isotrivial} \Longleftrightarrow
f_*\omega_{X_0/C_0} \mbox{ is a unitary flat line bundle on } C_0.$

\item Since $f$ is semi-stable, we have
$$ f_*\omega_{X_0/C_0} \mbox{ is unitary flat on } C_0
\Longleftrightarrow f_*\omega_{X/C}=\omega^{-1}_{C} \mbox{ is
unitary flat on } C .$$

\item $\omega^{-1}_{C} \mbox{ is unitary flat on } C
\Longleftrightarrow 0=\deg\omega^{-1}_C \Longleftrightarrow C \mbox{
is elliptic}.$
\end{itemize}

\item$(1)\Longrightarrow(3)$ : Otherwise, $C=\P^1$ since $C$ is
elliptic will induce that $f$ is isotrivial by the above agruments.
If $\#S\le 2$ then the global monodromy is generated by the locally
mondromies around points in $S.$ Since Deligne's complete reducible
theorem(cf.\cite{Del}) says that the global monodromy is
semi-simple, it implies that all local monodromies are identity,and
so the global monodromy is trivial. Then the VHS
$R^{n-1}f_*(\Q_{X_0})$ is trivial,  which contradicts to that $f$ is
non-isotrivial by the infinitesimal Torelli theorem.

\item $(3)\Longrightarrow(1)$  is obvious since that $f$ is isotrivial is
equivalent to that $C$ is elliptic.

\end{myenumii}

\end{proof}
\begin{myrem}
The proof actually implies the following equivalent conditions:
\begin{myenumiii}
    \item  $f$ is an isotrivial family;
    \item  $f_*\omega_{X/C}$ is an unitary flat line bundle on $C;$
    \item  $C$ is an elliptic curve.
\end{myenumiii}
%Therefore, there is an interesting phenomena : Even $f$ is
%isotrivial it is impossible that $C=\P^1$ with $\#S \leqq 2.$
\end{myrem}

%\begin{proposition}Let $f:X \to Y$ be a semi-stable
%family between two projective manifolds smooth over $Y_0$ such
%that $X$ has trivial canonical bundle and $S=Y\setminus Y_0$ is a
%reduced normal crossing divisor. Let $\mathfrak{M}_K$ be the
%quasi-projective coarse moduli scheme for the set of polarized
%projective manifolds with certain type 'K' as the general fibre of
%$f.$

%Assume that the moduli morphism $\eta_f: Y_0 \to \mathfrak{M}_K$
%is generically finite. Then, for any generic smooth rational curve
%$C \subset Y,$ $\# C \cap S \geq 3.$
%\end{proposition}

%From the results of the analogue Shafarevich conjecture for
%Calabi-Yau, we know
%\begin{corollary}[cf.]
%If $\mathfrak{M}_K$ is the coarse moduli space of polarized
%Calabi-Yau, then from the {\it rigidity} criterion by Viehweg-Zuo
%and Liu-Todorov-Yau-Zuo, any local monodromy $\rho \in \pi_1(C
%\setminus S\cap C)$ around point in $S\cap C$ can not be maximal
%unipotent. Also the restricted subfamily over $C$ can not be a
%Lefschetz pencil .
%\end{corollary}
%\begin{proof} Because $Y$ is a {\it rationally connected} manifold, any generic smooth
%rational curve can move freely.
%\end{proof}

%\begin{observation}[Existence of non-rigid Calabi-Yau pencil]In the real (physical)
%world, there exist Calabi-Yau manifolds fibred by lower
%dimensional Calabi-Yau varieties over a smooth projective manifold
%such that the induced moduli maps are generically finite.
%Therefore, our example shows that there are many non rigid
%families of Calabi-Yau over $\P^1.$
%\end{observation}

\section{Dimension counting for fibered Calabi-Yau manifolds}% with semi-stable Fibre Structures}
\subsection*{Calabi-Yau manifolds fibred by Abelian varieties}
%\subsubsection*{Fibre space over $\P^1$}

We now use Theorem \ref{main-theorem-1} and its corollaries to deal
with Calabi-Yau manifolds over $\C\P^1$ fibred by semi-stable
Abelian varieties. We then have the dimension of such Calabi-Yau
manifold $\leq 3,$ see Theorem \ref{main-theorem-3}.

\begin{proof}[Proof of Theorem \ref{main-theorem-3}]
Consider the Higgs bundle $(E,\theta)$ corresponding to
$R^1f_*\Q_{X_0}.$ The semi-stability of $f$ shows that all local
monodromies of $R^1f_*\Q_{X_0}$ are unipotent, thus there is the
Deligne canonical extension
$$\overline{E}= f_*\Omega^1_{X/\P^1}(\log \Delta)\bigoplus
R^1f_*(\sO_X)$$ where  $\Delta=f^{*}(S),$ and both pieces are
locally free.
\begin{myenumi}
\item Let $n=\dim f^{-1}(t).$ Then, $n=\rk
f_*\Omega^1_{X/\P^1}(\log \Delta)$ as a general fibre is an
Abelian variety. It has been shown in \ref{kollar-flat-ample} and
\ref{application-of-kollar} that $f$ is a non-isotrivial family
and the bundle $f_*\Omega^1_{X/\P^1}(\log \Delta)$ is ample on
$\P^1.$ The Grothendieck splitting theorem says that there is a
decomposition:
$$f_*\Omega^1_{X/\P^1}(\log \Delta)=\bigoplus_{i=1}^n\sO_{\P^1}(d_i),$$
and all integers $d_i$ are positive by the ampleness of the
$f_*\Omega^1_{X/\P^1}(\log \Delta).$

On the other hand, the commutative diagram of morphisms
$$
\begin{CDS}
\wedge^n f_*\Omega^1_{X/\P^1}(\log \Delta) \> \neq 0 >>f_*(\wedge^n\Omega^1_{X/}(\log \Delta))\\
\V V = V \novarr \V V = V \\
\sO_{\P^1}(\sum_{i=1}^{n} d_i) \> \neq 0  >> f_*\omega_{X/\P^1}.
\end{CDS}
$$ induces that
$$n\leq\sum_{i=1}^{n} d_i\leq \deg
f_*\omega_{X/\P^1}.$$ %=\deg(f_*\omega_X\otimes
%\omega_{\P^1}^{-1})=\deg f_*\omega_X +2.$$
By Zariski main theorem, $f$ only having connected fibres is same
as $f_*\sO_X=\sO_{\P^1}.$ Hence $\deg f_*\omega_{X/\P^1}=2,$ and
so the dimension of a general fibre is less than $3.$

\item  If $X$ is a Calabi-Yau threefold then
$$f_*\Omega_{X/\P^1}^1(\log \Delta)=\sO_{\P^1}(1)\oplus
\sO_{\P^1}(1).$$ There is so called Arakelov-Yau inequality (cf.
\cite{Del8}\cite{VZ6}).
$$\deg f_*\Omega_{X/\P^1}^1(\log \Delta)\leq
\frac{\rk{f_*\Omega_{X/\P^1}^1(\log \Delta) }}{2}
\deg(\Omega^1_{\P^1}(\log S))=2g(\P^1)-2 +\#S.$$ Thus $\#S\geq 4,$
and if $\# S =4$ then there is an \'etale covering $\pi:Y'\to \P^1$
such that $f':X'=X\times_{\P^1} Y'\to Y'$ is isogenous over $Y'$ to
a product $ E\times_{Y'}\cdots\times_{Y'} E ,$ where $h: E\to Y'$ is
a family of semi-stable elliptic curves
reaching the Arakelov bound.%,i.e., the inequality becomes inequality
%for
%$$\deg h_*\omega_{E/Y'}\leq \frac{\rk{h_*\omega_{E/Y'}}}{2} \deg(\Omega^1_{Y'}(\log S))=\frac{2g(Y')-2+\#S}{2}.$$
%The arithmetic comes from the result of Viehweg-Zuo directly.

\item If $X$ is a $K3$ surface then $$f_*\Omega_{X/\P^1}^1(\log
\Delta)=f_*\omega_{X/\P^1}=\sO_{\P^1}(2).$$ We deduce $\# S \geq
6$ from the Arakelov-Yau inequality for weight one VHS:
$$\deg f_*\Omega_{X/\P^1}^1(\log \Delta)\leq \frac{\rk{f_*\Omega_{X/\P^1}^1(\log \Delta) }}{2} \deg(\Omega^1_{\P^1}(\log S))=\frac{\#S}{2}-1.$$

\item The modularity is due to recent results by Viehweg-Zuo (cf.
\cite{VZ5}).
\end{myenumi}
\end{proof}

Now we deal with any semi-stable fiberation $f:X\to Y$ with higher
dimension base $Y.$ We obtain the dimension of a general fibre is
bounded above by a constant  dependent on $Y.$
\begin{proof}[Proof of Theorem \ref{main-theorem-4}]
The proof follows from the next steps.
\begin{myenumi}
 \item By \ref{Kollar-decomposition}, $f_*\Omega^1_{X/Y}(\log \Delta)$  has no flat quotient.
 Moreover, $f_*\Omega^1_{X/Y}(\log \Delta)$ is ample on
 any {\it sufficiently general} curve in $Y.$ We always have:
$$\wedge^n f_*\Omega^1_{X/Y}(\log \Delta)=\det
f_*\Omega^1_{X/Y}(\log \Delta),$$ where $n$ is the dimension of a
general fibre and also is the rank of the locally free sheaf
$f_*\Omega^1_{X/Y}(\log \Delta).$ On the other hand, the non-zero
map
$$ \wedge^n f_*\Omega^1_{X/Y}(\log \Delta) \> \neq 0
>>f_*(\wedge^n\Omega^1_{X/Y}(\log
 \Delta))=f_*\omega_{X/Y},$$
induces that the line bundle $\omega^{-1}_Y=f_*\omega^1_{X/Y}$ is
{\it big} and {\it nef}. Thus $Y$ is a rationally connected
projective manifold and so $\sU=0$ by \ref{main-theorem-1}.

\item By Proposition \ref{Kollar-fundamental-grp}, we have a very free
morphism $g:\P^1 \to \overline{M}$ which is sufficiently general,
i.e., the image curve $C:=g(\P^1)$ satisfies:
\begin{myenumii}
    \item $C$ intersects $S$ transversely;
    \item $\pi_1(C_0) \onto \pi_1(Y_0)\to 0$ is surjective where $C_0=C\cap
Y_0.$% and $D_\infty=C\cap S.$
\end{myenumii}
\item If $\dim Y \geq 3,$  $g$ is an embedding and we have a very
free smooth rational curve $C \subset Y.$ Since
$f_*\Omega^1_{X/Y}(\log \Delta)$ is ample over $C,$
$$f_*\Omega^1_{X/Y}(\log \Delta)|_C \cong \bigoplus_{i=1}^n\sO_{\P^1}(d_i)
\mbox{ with }  \forall d_i>0.$$ Denote by $l=-\omega_Y \cdot C.$ The
commutative diagram of morphisms
$$
\begin{CDS}
\wedge^n f_*\Omega^1_{X/Y}(\log \Delta)|_C \> \neq 0 >>f_*(\wedge^n\Omega^1_{X/Y}(\log \Delta))|_C\\
\V V = V \novarr \V V = V \\
\sO_{\P^1}(\sum_{i=1}^{n} d_i) \> \neq 0  >> f_*\omega_{X/Y}|_C.
\end{CDS}
$$
induces that
$$n\leq \sum_{i=1}^{n} d_i \leq \deg_C
f_*\omega_{X/Y}=-\deg_C(\omega_Y)=l.$$

\item If $\dim Y=1,$ $Y$ is then $\P^1$ and $l=2.$ If $\dim Y=2,$
then $Y$ is a smooth Del Pezzo surface and $g$ is an immersion by
the theorem \ref{Kollar-fundamental-grp}.
\begin{myenumii}
\item Let $\Gamma$ be the graph of the morphism $g,$ i.e.,
$$\Gamma=\{(x,y)\in \P^1\times X \ | \  y=g(x) \} \subset \P^1\times X.$$
$\Gamma$ is a smooth curve, actually it is isomorphic to $\P^1.$
We have
$$
\begin{TriCDA}
{\P^1\times X}  {\SW \mathrm{pr}_1 WW}{\SE E \mathrm{pr}_2 E}
{\P^1} { } {X}
\end{TriCDA}
$$
and the projections $\mathrm{pr}_1$ and $\mathrm{pr}_2$ both are
proper morphisms. Hence $\mathrm{pr}_1 : \Gamma \to g(\P^1)$ is a
finite morphism. Denote by $C:=g(\P^1).$ We then have a finite set
$B\subset C$ such that $g^{-1}(B)$ is a finite set in $\P^1$ and
$$g: \P^1\setminus g^{-1}(B) \longrightarrow C\setminus B$$ is an
\'etale covering.

\item  Since the real-codimension of $B$ in $Y$ is four, we have a
natural isomorphism of topological fundamental groups
$\pi_1(Y_0\setminus B)\> \cong >> \pi_1(Y_0),$ and so a surjective
homomorphism $\pi_1(C_0\setminus B) \onto \pi_1(Y_0) \to 0.$ Denote
by $$T_0:=\P^1-g^{-1}(B\cup (C-C_0))\mbox{ and }\phi:=g|_{T_0}.$$ We
have an \'etale covering $\phi: T_0 \to C_0\setminus B,$ and an
injective
$$0\>>> \sO_{C_0\setminus B} \>>>
\phi_*\sO_{T_0}.$$ Without lost generality we may assume that $\phi$
is a Galois covering, then the above short sequence has a split and
so
$$\phi_*\sO_{T_0}=\sO_{C_0\setminus B}\oplus \mbox{ Galois
conjugates}.$$

 \item Let $\sF=f_*\Omega^1_{X/Y}(\log \Delta)$ and  $\sL:=g^*(\sF|_C).$
Then,
$$\sL=\bigoplus_{i=1}^n \sO_{\P^1}(d_i) \mbox{ with } d_i \geq 0 \ \forall i.$$
We claim that all $d_i>0.$

Suppose that one $d_i=0,$ then $\phi_*(\sL|_{T_0})$ has a nonzero
flat quotient $$\phi_*(\sL|_{T_0}) \onto \sO_{C_0\setminus B} \to
0.$$ Thus, $\sF|_{C_0 \setminus B}$ has a nonzero flat quotient
since we have
$$0\>>> \sF|_{C_0 \setminus B} \>>> \phi_*\phi^*( \sF|_{C_0 \setminus B})=\phi_*(\sL|_{T_0}).$$
and the projection formula
$$ \phi_*\phi^*( \sF|_{C_0 \setminus B})= \sF|_{C_0
\setminus B}\otimes \phi_*(\sO_{T_0}).$$

The subjectivity of $\pi_1(C_0\setminus B) \onto \pi_1(M) \to 0$
implies that $\sF|_M$ has a flat quotient, and so
$f_*\Omega^1_{X/Y}(\log \Delta)$ itself has a unitary flat quotient.
It is a contradiction. Hence, all integers $d_i$ are positive and
$$n\leq -\deg_{\P^1}g^*\omega_Y=-g_*[\P^1]\cdot \omega_Y=l.$$
\end{myenumii}
\end{myenumi}
\end{proof}
%\noindent{\bf Remark.} Since $g$ is a very free morphism, $l\geq
%\dim Y+1.$ However, it seems that we can not find a very free
%rational curve with $l=\dim Y +1$ in the theorem.
%$$g^*T_Y=\bigoplus_{i=1}^{\dim Y} \sO_{\P^1}(a_i) \mbox{ with } a_1\geq a_2 \geq \cdots \geq a_{\dim Y}\geq 1
%\mbox{ and } a_1\geq2. $$

\subsection*{Calabi-Yau manifolds fibred by hyperk\"ahler varieties}Using similar methods as the case
of Calabi-Yau  manifolds fibred by Abelian varieties, we have the following result : % (\ref{dim-hyperkahler} and \ref{main-theorem-3})
\begin{theorem}\label{dim-hyperkahler} Let $f: X\to \P^1$ be a semi-stable family fibred by
hyperk\"ahler varieties with $$f:X_0=f^{-1}(C_0) \to C_0$$ smooth,
$S=\P^1 \setminus C_0,$ and $\Delta=f^{-1}(S).$ Assume that $X$ is a
projective manifold with  trivial canonical sheaf $\omega_X$ and
$H^2(X,\sO_X)=0.$ Then, $f$ is nonisotrivial with $\#S \geq 3$ and
the dimension of a general fibre is four.
\end{theorem}

\begin{proof}
Consider the Higgs bundle $(E,\theta)$ corresponding to
$R^2f_*\Q_{X_0}.$ The semi-stability of $f$ shows that there is
the Deligne canonical extension:
$$\overline{E}= f_*\Omega^2_{X/\P^1}(\log \Delta)\bigoplus
R^1f_*\Omega^1_{X/\P^1}(\log \Delta)\bigoplus R^2f_*(\sO_X),$$
such that $f_*\Omega^2_{X/\P^1}(\log \Delta), R^2f_*(\sO_X)$ are
line bundles. % since $\rank f_*\Omega^2_{X/\P^1}(\log \Delta)=1.$
By \ref{kollar-flat-ample},
$$f_*\Omega^2_{X/\P^1}(\log \Delta)=\sO_{\P^1}(d)$$ and $d$ is a
positive integer.

Let  $F$ be a general fibre and let $ n=\dim F/2.$ The commutative
diagram
$$
\begin{CDS}
 \mathrm{Sym}^n f_*\Omega^2_{X/\P^1}(\log \Delta) \> \neq 0 >>f_*(\wedge^n\Omega^1_{X/\P^1}(\log \Delta)) \\
\V V = V \novarr \V V = V \\
\sO_{\P^1}(nd) \> \neq 0  >>  f_*\omega_{X/\P^1},
\end{CDS}
$$
%and the non-zero morphism $$\sO_{\P^1}(nd) \> \neq 0 >> f_*\omega_{X/\P^1}=\sO_{\P^1}(2)$$
induces $n \leq  \deg f_*\omega_{X/\P^1}=2.$ Thus $\dim F=4$ by
the definition, and so
$$f_*\Omega^2_{X/\P^1}(\log \Delta)= \sO_{\P^1}(1).$$

%For the VHS $R^4f_*(\C_{X_0}),$ the Arakelov-Yau inequality  and
%the infinitesimal Torelli theorem shows that
%$$\deg
%f_*\omega_{X/\P^1}\le \frac{\dim
%F}{2}\rank(f_*\omega_{X/\P^1})\deg \Omega_{\P^1}^1(\log S)=2(\#S
%-2),$$ and so $\#S \geq 3.$
$\#S\ge 3$ is a well-known result since $f$ is nonisotrivial.
\end{proof}
%\noindent{\bf Remark.} Here, $\overline{\theta}^{2,0}$ has to be
%an injection.

%Use the similar methods as in \ref{dim-hyperkahler} and \ref{main-theorem-3}, we have:

\begin{theorem}\label{main-theorem-5}
Let $f:X\to Y$ be a semi-stable family of hyperk\"ahler varieties
between two non-singular projective varieties with
$f:X_0=f^{-1}(Y_0) \to Y_0$ smooth, $S=Y \setminus Y_0$ a reduced
normal crossing divisor and $\Delta=f^{*}(S)$ a relative reduced
normal crossing divisor in $X.$ Assume that the period map for the
VHS $R^2f_*(\Q_{X_0})$
is injective at one point in $Y_0$ and %ing induced moduli morphism is generally finite and Assume that
$X$ has trivial canonical sheaf $\omega_X$ with
$H^0(X,\Omega^2_X)=0.$ Then, the dimension of a general fibre is
bounded above by a constant depending on $Y.$
\end{theorem}
\begin{remark}
If a fibration $f: X\to Y$ has $\dim X=2n$ and $X$ is a projective
irreducible sympletic manifold, then $\dim Y=n$(cf. \cite{Matsu}).
Here we obtain some necessary conditions for a question asked by
Naichung Conan Leung in CUHK algebraic geometry working seminar :
{\it Does there exist a Calabi-Yau manifold fibred by Abelian
varieties (resp. hyperk\"ahler varieties)? }
%We believe that the dimension of a general fibre should be bounded above by a constant which only depends  on $\dim Y.$ \\
\end{remark}

\centerline{{\bf Acknowledgements}}

\vspace{0.5cm} We wish to thank Yi Hu for telling us the results in
\cite{Matsu}, and thank Eckart Viehweg for valuable suggestions on
moduli spaces and the positivity of relative dualizing sheaves. The
first author is very grateful to Shing-Tung Yau for his constant
encouragement.

\bibliographystyle{plain}

\end{document}